\documentstyle[12pt]{article}

\def\thebibliography#1{\section*{References}\list
{$^{\arabic{enumi}}$}{\settowidth\labelwidth{#1}\leftmargin\labelwidth
\advance\leftmargin\labelsep
\usecounter{enumi}}
\def\newblock{\hskip .11em plus .33em minus .07em}
\sloppy\clubpenalty4000\widowpenalty4000
\sfcode`\.=1000\relax}

\def\op#1{\mathop{\fam0 #1}\limits}

\newcommand{\pr}{{\rm pr}}
\newcommand{\Ker}{{\rm Ker\,}}
\newcommand{\im}{{\rm Im\, }}

\newcommand{\beq}{\begin{equation}}
\newcommand{\eeq}{\end{equation}}
\newcommand{\ben}{\begin{eqnarray}}
\newcommand{\een}{\end{eqnarray}}
\newcommand{\be}{\begin{eqnarray*}}
\newcommand{\ee}{\end{eqnarray*}}
\newcommand{\bea}{\begin{eqalph}}
\newcommand{\eea}{\end{eqalph}}

\newcommand{\cT}{{\cal T}}
\newcommand{\cL}{{\cal L}}
\newcommand{\cE}{{\cal E}}
\newcommand{\cH}{{\cal H}}
\newcommand{\bL}{{\bf L}}
\newcommand{\R}{{\bf R}}
\newcommand{\al}{\alpha}

\newcommand{\dl}{\delta}
\newcommand{\la}{\lambda}

\newcommand{\f}{\phi}

\newcommand{\Om}{\Omega}
\newcommand{\m}{\mu}

\newcommand{\g}{\gamma}
\newcommand{\G}{\Gamma}

\newcommand{\si}{\sigma}
\newcommand{\vt}{\vartheta}
\newcommand{\bom}{{\bf\Omega}}
\newcommand{\w}{\wedge}
\newcommand{\wt}{\widetilde}
\newcommand{\wh}{\widehat}
\newcommand{\ol}{\overline}
\newcommand{\dr}{\partial}

\newcommand{\ap}{\approx}

\let\ssection=\section
\renewcommand{\section}{\setcounter{equation}{0}\ssection}

\newcounter{eqalph}
\newcounter{equationa}
\newcounter{theorem}
\newcounter{proposition}
\newcounter{lemma}
\newcounter{corollary}
\newcounter{definition}

\newenvironment{eqalph}{\stepcounter{equation}
\setcounter{equationa}{\value{equation}}
\setcounter{equation}{0}

\begin{eqnarray}}{\end{eqnarray}\setcounter{equation}{\value{equationa}}}

\def\thedefinition{\arabic{definition}}

\newenvironment{proof}{\noindent 
{\it Proof:}}{\medskip}
\newenvironment{rem}{\medskip\noindent{\it Remark:}}{\medskip}

\newenvironment{prop}{\refstepcounter{definition} 
\bigskip\noindent{\it Proposition \thedefinition:}}{\medskip}
\newenvironment{lem}{\refstepcounter{definition} 
\bigskip\noindent{\it Lemma \thedefinition:}}{\medskip}

\newenvironment{defi}{\refstepcounter{definition} 
\bigskip\noindent{\it Definition \thedefinition:}}{\medskip} 

\begin{document}
\hbox{}
 
\centerline{\large\bf NON-SYMPLECTIC GEOMETRY OF}
\bigskip

\centerline{\large\bf FIRST ORDER TIME-DEPENDENT MECHANICS}
\bigskip
 
\centerline{\sc Gennadi A. Sardanashvily}
\medskip
 
\centerline{Department of Theoretical Physics, Physics Faculty,}
 
\centerline{Moscow State University, 117234 Moscow, Russia}

\centerline{E-mail: sard@grav.phys.msu.su}
\bigskip
 
\begin{abstract}
The usual formulation of time-dependent mechanics implies a given splitting
$Y=R\times M$ of an event space $Y$. This splitting, however, is broken by
any time-dependent transformation, including transformations between inertial
frames. The goal is the frame-covariant formulation of time-dependent 
mechanics on a bundle $Y\to R$ whose fibration $Y\to M$ is not fixed. Its 
phase space is the vertical cotangent bundle $V^*Y$ provided with the 
canonical 3-form and the corresponding canonical Poisson structure. An event 
space of relativistic mechanics is a manifold $Y$ whose fibration $Y\to R$ 
is not fixed.
\end{abstract}
\bigskip
\bigskip

\noindent 
{\bf I. INTRODUCTION}
\bigskip

Symplectic technique is well-known to provide adequate mathematical
formulation of autonomous mechanics, when Hamiltonians are independent
of time.$^{1-4}$ 
Its canonical example is a mechanical system whose phase space
is the cotangent bundle $T^*M$ of an event manifold $M$. 
This phase space is provided with the canonical symplectic form
\be
\Om= dp_i\w dy^i, 
\ee
written with respect to the holonomic coordinates $(y^i, p_i=\dot q_i)$ on 
$T^*M$. A Hamiltonian $\cH$ is 
defined as a real function on $T^*M$. The motion trajectories
are integral curves of the Hamiltonian vector field 
$\vt=\vt_i\dr^i +\vt^i\dr_i$ on $T^*M$ which obeys the Hamilton equations
\ben
&&\vt\rfloor\Om=d\cH, \label{m2}\\
&& \vt^i=\dr^i\cH, \qquad \vt_i=-\dr_i\cH. \nonumber
\een

The usual formulation of time-dependent mechanics 
implies a splitting
$Y=\R\times M$
of the event manifold $Y$ and the corresponding
splitting
$\R\times T^*M$ of the phase 
space.$^{5-10}$ The phase space is provided with the
pull-back $\pr^*_2\Om$ of the symplectic form on $T^*M$.
By a time-dependent Hamiltonian is meant a real function on 
$\R\times T^*M$, 
while motion trajectories
are integral curves of the time-dependent Hamiltonian vector field 
\be
\vt: \R\times T^*M\to TT^*M
\ee
which obeys the Hamilton equations (\ref{m2}).
The problem is that the above-mentioned splittings are broken by any
time-dependent transformation, including transformations of inertial frames.
Therefore, the form $\pr^*_2\Om$ on the phase space of time-dependent
mechanics fails to be canonical.$^{11}$

We will formulate first order time-dependent mechanics as a
particular field theory, when an event space is a fibred
manifold $Y\to\R$, coordinatized by $(t,y^i)$.$^{5,12-14}$
 The configuration space is the
first order jet manifold $J^1Y$ of sections of $Y\to\R$, which is
provided with the adapted coordinates $(t,y^i,y^i_t)$. 
There is the canonical monomorphism 
\be
\la: J^1Y\hookrightarrow TY, \qquad \la=\dr_t +y^i_t\dr_i,
\ee
over $Y$. It is easy to see that $\pi^1_0:J^1Y\to Y$ is an affine bundle
modelled over the vertical tangent bundle $VY\to Y$.
For the sake of simplicity, we will 
identify $J^1Y$ with the corresponding subbundle of $TY$. 

The 1-dimensional reduction of 
polysymplectic Hamiltonian formalism$^{13-16}$
provides the adequate mathematical
formulation of time-dependent Hamiltonian mechanics on 
the Legendre bundle $\pi_\Pi:\Pi= V^*Y\to Y$. 
The phase space $V^*Y$ is endowed with the canonical 3-form
\beq
\bom=dp_i\w dy^i\w dt, \label{m66}
\eeq
written with respect to the holonomic coordinates $(t,y^i,p_i=\dot
y_i)$ on $V^*Y$. 

\begin{rem}
Unless otherwise stated, the base $\R$ is parameterized by the
coordinates $t$ with transition functions
$t'=t+$const. Relative to these coordinates, $\R$ is equipped with the
 standard vector field $\dr_t$
and the standard 1-form
$dt$, which is  also the volume element on $\R$. This is not the
case of relativistic mechanics.
\end{rem}

The following peculiarities of  
time-dependent Hamiltonian  mechanics should be emphasized.

(i) The form $\bom$ (\ref{m66}) defines the canonical degenerate
Poisson structure on the phase space $V^*Y$.

(ii) A Hamiltonian is not a function on a
phase space. As a consequence, the evolution equation is not reduced to a
Poisson bracket, and integrals of motion cannot be defined as
functions in involution with a Hamiltonian. 

(iii)  Hamiltonian and Lagrangian
formulations of time-dependent mechanics are equivalent in the case of
hyperregular Lagrangians. 
A degenerate Lagrangian requires a set of associated Hamiltonians and
Hamilton equations 
in order to exhaust all solutions of the Lagrange equations. 

(iv) A connection
\beq
\G:Y\to J^1Y\subset TY, \qquad \G=\dr_t +\G^i\dr_i,\label{1005}
\eeq
on the event space $Y\to\R$ defines a reference frame.
There is one-to-one correspondence between the connections $\G$ 
and atlases of constant local
trivializations of $Y\to\R$ such that the transition functions
$y^i\to y'^i(y^j)$ are independent of $t$.$^{16,17}$
Thus, we obtain the frame-covariant formulation of time-dependent
mechanics, that enables us to describe phenomena related to different
reference frames.

For the sake of simplicity, $Y\to\R$ is assumed to be a bundle with
a typical fibre $M$. It is trivial. Different trivializations
\beq
Y\cong \R\times M \label{m33}
\eeq
differ from each other in fibrations $Y\to M$, while the fibration
$\pi:Y\to\R$ is once for all.
Given a trivialization (\ref{m33}),
there are the corresponding splittings of the configuration and phase spaces
\be
J^1Y\cong \R\times TM, \qquad V^*Y\cong \R\times T^*M. 
\ee

If a fibration  $Y\to \R$ of an event space $Y$ is not fixed, we
obtain the general formulation of relativistic mechanics, including 
Special Relativity on $Y=\R^4$.$^{16}$ 
Its configuration space is the first order jet manifold
$J^1_1Y$ of 1-dimensional submanifolds of $Y$. 

All manifolds throughout are assumed to be paracompact and connected.
\bigskip

\noindent 
{\bf 2. CANONICAL POISSON STRUCTURE}
\bigskip

The Legendre
bundle $V^*Y$ of time-dependent mechanics is provided with the canonical
Poisson structure as follows.             
Let us consider the cotangen bundle $T^*Y$ 
with the holonomic 
coordinates $(t,y^i,p_i,p)$, which is the homogeneous Legendre bundle
of time-dependent mechanics. It admits the
canonical Liouville form
\beq
\Xi=pdt +p_idy^i \label{m91}
\eeq
and the canonical symplectic form
\be
d\Xi=dp\w dt +dp_i\w dy^i. 
\ee
The corresponding Poisson bracket on the space $C^\infty(T^*Y)$ of
functions on $T^*Y$ reads
\beq
\{f,g\} =\dr^pf\dr_tg - \dr^pg\dr_tf +\dr^if\dr_ig-\dr^ig\dr_if. \label{m116}
\eeq

Let us
consider the subspace of
$C^\infty(T^*Y)$ which comprises the pull-backs of functions on $V^*Y$ by the
projection
$T^*Y\to V^*Y$. This subspace is closed under the
Poisson bracket (\ref{m116}). Then there
exists the canonical Poisson structure
\beq
\{f,g\}_V = \dr^if\dr_ig-\dr^ig\dr_if \label{m72}
\eeq
on $V^*Y$ induced by (\ref{m116}).$^3$ 
The corresponding Poisson bivector 
\be
w(df,dg)= \{f,g\}_V 
\ee
on  $V^*Y$ is vertical with respect to 
the fibration $V^*Y\to \R$, and reads
\be
w^{ij} =0, \qquad w_{ij}=0, \qquad w^i{}_j=1. 
\ee
A glance at this expression shows that the holonomic coordinates on $V^*Y$ 
are canonical for the Poisson structure (\ref{m72}), which is
regular and degenerate. 

Given the Poisson bracket (\ref{m72}), the Hamiltonian vector field $\vt_f$ 
of a function
$f$ on $V^*Y$, defined by the relation $\{f,g\}_V=\vt_f\rfloor dg$,
$\forall g\in C^\infty(V^*Y)$, is the vertical vector field 
\beq
\vt_f = \dr^if\dr_i- \dr_if\dr^i \label{m73}
\eeq
on $V^*Y\to \R$. Hence, the characteristic distribution of the
Poisson structure (\ref{m72}), generated by Hamiltonian vector fields, 
is precisely the vertical tangent bundle $VV^*Y$ of
$V^*Y\to \R$. 

By virtue of the well-known theorem,$^4$
the Poisson structure (\ref{m72}) defines the symplectic foliation on
$V^*Y$ which coincides with the fibration
$V^*Y\to \R$. 
The symplectic forms on the fibres
of
$V^*Y\to\R$ are the pull-backs  
\be
\Om_t=dp_i\w dy^i
\ee
of the canonical symplectic form on the typical fibre
$T^*M$ of $V^*Y\to \R$ with respect to trivialization morphisms.$^{11}$

The Poisson structure (\ref{m72}) can be introduced in a different
way. The Legendre bundle $V^*Y$ admits the 
canonical closed 3-form (\ref{m66}),
which is the particular case of the polysymplectic
form.$^{14-16}$ Then every function
$f$ on $V^*Y$ defines the corresponding Hamiltonian vector field 
$\vt_f$ (\ref{m73}) by the relation
\be
\vt_f\rfloor\bom = df\w dt,
\ee
while the Poisson bracket (\ref{m72}) is recovered by the condition
\be
\{f,g\}_Vdt=\vt_g\rfloor\vt_f\rfloor\bom.
\ee 
\bigskip

\noindent 
{\bf 3. HAMILTONIAN FORMS}
\bigskip

Following general polysymplectic formalism,
we say that a connection
\beq
\g=\dr_t + \g^i\dr_i +\g_i\dr^i \label{m115}
\eeq
on $V^*Y\to \R$ is locally Hamiltonian if the exterior form
$\g\rfloor\bom$ is closed, i.e.,
\beq
\bL_\g\bom=d(\g\rfloor\bom)=0 \label{m68}
\eeq
where $\bL$ denotes the Lie derivative.

For instance, every connection
$\G$ on the bundle
$Y\to \R$ gives rise to the locally Hamiltonian connection 
\be
\wt\G =\dr_t +\G^i\dr_i -p_i\dr_j\G^i \dr^j
\ee
such that
\beq
\wt\G\rfloor\bom=dH_\G, \qquad H_\G=p_idy^i -p_i\G^idt. \label{m61}
\eeq

Locally Hamiltonian connections constitute an affine space modelled over 
the linear space of 
vertical vector fields $\vt$ on $V^*Y\to \R$
which obey the same condition (\ref{m68}), and are locally Hamiltonian vector 
fields as follows.

\begin{lem}\label{exact}
Every
closed form $\g\rfloor\bom$ on $V^*Y\to \R$ is exact. 
\end{lem}

\begin{proof}
Let us consider the decomposition 
\beq
\g=\wt\G +\vt \label{b4211}
\eeq
where $\G$ is a connection on $Y\to \R$, while $\vt$ satisfies the relation
$d(\vt\rfloor\bom)=0$.
It is easily seen that $\vt\rfloor\bom=\si\w dt$ 
where  $\si$ is a 1-form.  Using 
properties of the De Rham
cohomology groups of a manifold product, one can show that 
every closed 2-form $\si\w dt$ on $V^*Y$ is
exact, and so is $\g\rfloor\bom$. 
Moreover, in 
accordance with the relative Poincar\'e lemma, we can write locally
$\vt\rfloor\bom =df\w dt$.
\end{proof}
 
\begin{defi}
A 1-form $H$ on the Legendre bundle
$V^*Y$ is called locally Hamiltonian if
\be
\g\rfloor\bom=dH
\ee
for a connection
$\g$ on $V^*Y\to\R$.
\end{defi}

By virtue of
Proposition \ref{exact}, there is 
one-to-one correspondence between the locally Hamiltonian connections and the
locally Hamiltonian forms considered throughout modulo closed forms.

\begin{defi} By a Hamiltonian form $H$ on
the Legendre bundle
$V^*Y$ is meant the pull-back
\beq
H=h^*\Xi= p_i dy^i -\cH dt \label{b4210}
\eeq
of the Liouville form $\Xi$ (\ref{m91}) on $T^*Y$ by a section $h$ of the
bundle $T^*Y\to V^*Y$. 
\end{defi}

Any connection $\G$ on $Y\to \R$ defines the
Hamiltonian form $H_\G$ (\ref{m61}) on $V^*Y$, and every Hamiltonian
form on $V^*Y$ admits the splitting
\beq
H =p_idy^i -(p_i\G^i +\wt{\cH}_\G) dt  \label{m46}
\eeq
where $\G$ is a connection on $ Y\to\R$ and $\wt\cH_\G$ is a real 
function on $V^*Y$. 

Given a trivialization (\ref{m33}) of $Y\to \R$ 
the Hamiltonian form (\ref{b4210}) looks like the well-known 
Poincar\'e--Cartan integral.$^2$  However, the Hamiltonian $\cH$ in
the expression (\ref{b4210}) is not a function. A glance at the
splitting (\ref{m46}) shows that Hamiltonians (and Hamiltonian forms) 
constitute an affine space modelled over the linear space of functions on
$V^*Y$. 

\begin{prop}\label{genform}
Locally Hamiltonian forms are Hamiltonian forms locally.
\end{prop}

\begin{proof}
Given locally Hamiltonian forms $H_\g$ and $H_{\g'}$, their difference 
\be
\si=H_\g -H_{\g'}, \qquad d\si= (\g-\g')\rfloor\bom,
\ee
is a 1-form on $V^*Y$ such that the 2-form $\si\w dt$ is closed and,
consequently, exact.
In accordance with the relative Poincar\'e lemma, this condition implies
that $\si=fdt + dg$
where $f$ and $g$ are local functions on $V^*Y$.
Then it follows from the splitting (\ref{b4211})
that, in a neighbourhood of every point $p\in V^*Y$, a locally Hamiltonian 
form $H_\g$
coincides with the pull-back of the Liouville form $\Xi$ on $T^*Y$ 
by the local section
\be
(t,y^i,p_i)\mapsto (t,y^i,p_i,p=-p_i\G^i+f)
\ee
of $T^*Y\to V^*Y$ where $f$ is a local function on $V^*Y$. 
\end{proof}

\begin{prop} Conversely, let $H$ be 
a Hamiltonian form $H$ on the Legendre bundle $V^*Y$.
There exists a unique connection $\g_H$ on $V^*Y\to \R$, called the
Hamiltonian connection, such
that $\g_H\rfloor\bom= dH$. 
\end{prop}

\begin{proof}
As in the polysymplectic case, let us introduce the 
Hamilton operator on the phase space $V^*Y$. This reads
\be
\cE_H=dH-\wh\bom=[(y^i_t-\dr^i\cH) dp_i
-(p_{ti}+\dr_i\cH) dy^i]\w dt
\ee
where $\wh\bom$
is the pull-back of the canonical form $\bom$ (\ref{m66}) onto $J^1V^*Y$. 
It is readily observed that the kernel of $\cE_H$ is an affine
subbundle of the Legendre bundle $V^*Y\to Y$. Therefore, its global section
\beq
\g_H=\dr_t +\dr^i\cH\dr_i -\dr_i\cH\dr^i \label{m57}
\eeq
always exists. This is the unique solution of the
first order differential Hamilton equations
\beq
y^i_t =\dr^i\cH, \qquad  p_{ti} =-\dr_i\cH \label{m41}
\eeq
on $V^*Y$, and is a Hamiltonian connection for the Hamiltonian form $H$.
\end{proof}

The integral curves of the Hamiltonian
connection (\ref{m57}) are
classical solutions of the Hamilton equations (\ref{m41}). 
Conversely, since the bundle  $\Ker\cE_H\to Y$ is affine,
every classical solution $r:\R\to V^*Y$
of the Hamilton equations (\ref{m41})
can be extended to a Hamiltonian connection for $H$.

Hamiltonian connections $\g_H$ (\ref{m57}) form an
affine space modelled over the linear space of Hamiltonian vector fields
(\ref{m73}).

\begin{rem}
Note that the Hamilton equations (\ref{m41}) can be
introduced without appealing to the Hamilton operator. They are equivalent
to the relation
\be
r^*(u\rfloor dH)=0
\ee
which is assumed to hold for any vertical vector field $u$ on $V^*Y\to\R$.
\end{rem}

With a Hamiltonian form $H$ (\ref{m46}) and the corresponding 
Hamiltonian connection $\g_H$ (\ref{m57}), 
we have the Hamilton evolution equation
\beq
d_{Ht}f=g_H\rfloor df=(\dr_t +\dr^i\cH\dr_i -\dr_i\cH\dr^i)f
\label{m59}
\eeq
on functions on the Legendre bundle $V^*Y$. Substituting a classical solution
of the Hamilton equations (\ref{m41}) in (\ref{m59}), we
obtain the time evolution of the function $f$.
Given the splitting (\ref{m46}) of a Hamiltonian form $H$, the
Hamilton evolution equation (\ref{m59}) is
brought into the form 
\beq
d_{Ht}f= \dr_tf +(\G^i\dr_i -\dr_i\G^jp_j\dr^i)f
+\{\wt\cH_\G,f\}_V. \label{m96}
\eeq
A glance at this expression shows that the Hamilton evolution
equation in time-dependent mechanics does not reduce to the Poisson
bracket. This fact may
be relevant to the quantization problem. The second term in the
right-hand  side of the equation (\ref{m96}) remains classical. 
\bigskip

\noindent 
{\bf 4. CANONICAL TRANSFORMATIONS}
\bigskip

Canonical transformations in time-dependent mechanics are not
compatible with the fibration $V^*Y\to Y$. 

\begin{defi}\label{can3} By a
canonical automorphism is meant an
automorphism 
$\rho$ over $\R$ of the bundle
$V^*Y\to \R$ which preserves the canonical 
Poisson structure
(\ref{m72}) on $V^*Y$, i.e., 
\be
\{f\circ\rho,g\circ\rho\}_V=\{f,g\}_V\circ\rho
\ee
and, equivalently, the canonical form 
$\bom$ (\ref{m66}) on $V^*Y$, i.e., $\bom=\rho^*\bom$. 
\end{defi}

The bundle
coordinates on $V^*Y\to \R$ are called canonical
if they are canonical for the Poisson structure (\ref{m72}). 
Canonical coordinate transformations satisfy the relations
\be
\frac{\dr {p'}_i}{\dr p_j}\frac{\dr {y'}^i}{\dr p_k}
-\frac{\dr {p'}_i}{\dr p_k}\frac{\dr {y'}^i}{\dr p_j}=0, \quad
\frac{\dr {p'}_i}{\dr y^j}\frac{\dr {y'}^i}{\dr y^k}
-\frac{\dr {p'}_i}{\dr y^k}\frac{\dr {y'}^i}{\dr y^j}=0, \quad
\frac{\dr {p'}_i}{\dr p_j}\frac{\dr {y'}^i}{\dr y^k}
-\frac{\dr {p'}_i}{\dr y^j}\frac{\dr {y'}^i}{\dr p_k} = \dl^k_j. 
\ee
By definition, the
holonomic coordinates on $V^*Y$ are the canonical ones.
Accordingly,  holonomic automorphisms 
\be
(y^i,p_i)\mapsto (y'^i, p'_i=\frac{\dr y^j}{\dr y'^i}p_j)
\ee
of the Legendre bundle $V^*Y\to Y$
induced by the vertical automorphisms of $Y\to\R$ are also canonical.

\begin{prop}\label{genham}
Canonical automorphisms send locally Hamiltonian connections onto the
locally Hamiltonian
ones (and, consequently, locally Hamiltonian forms onto each other). 
\end{prop}

\begin{proof}
If $\g$ is a locally Hamiltonian connection for $H$, we have 
\be
T\rho(\g)\rfloor\bom=(\rho^{-1})^*(\g\rfloor\bom) =d((\rho^{-1})^*H).
\ee
\end{proof}

\begin{prop}\label{p13.5} Let $\g$ be a complete 
locally Hamiltonian connection on $V^*Y\to \R$, i.e., the vector
field (\ref{m115}) is complete. 
There exist canonical coordinate
transformations which bring all components of $\g$ to zero, i.e.,
$\g=\dr_t$.
\end{prop}

\begin{proof}
A glance at the relation (\ref{m68}) shows that each locally Hamiltonian
connection $\g$
is the generator of a local 1-parameter group 
$G_\g$ of canonical automorphisms of
$V^*Y\to \R$. 
Let $V^*_0Y$ be the fibre of $V^*Y\to \R$ at $0\in \R$. 
Then  canonical coordinates of $V^*_0Y$ dragged along integral
curves of the complete vector field
$\g$ satisfy the statement of the proposition. 
\end{proof}

In particular, let $H$ be a Hamiltonian form (\ref{m46}) such that the
corresponding Hamiltonian connection $\g_H$ (\ref{m57}) is
complete. By virtue of Proposition \ref{p13.5}, there exist canonical
coordinate transformations which bring the Hamiltonian $\cH$ into zero. Then 
the corresponding Hamilton equations reduce to the equilibrium ones
\be
y^i_t=0, \qquad p_{ti}=0. 
\ee
Accordingly, any Hamiltonian form $H$ can be locally brought
into the form where $\cH=0$ by local canonical coordinate
transformations.
It should be emphasized that, in general, 
canonical automorphisms do not send Hamiltonian
forms onto Hamiltonian forms, but only locally.

Let $H$ be a Hamiltonian form (\ref{b4210}) on $V^*Y$. Given a canonical
automorphism $\rho$, we have
\be
d(\rho^*H -H)=0.
\ee
It follows that
\be
\rho^*H -H =dS
\ee
where $S$ is a local function on $V^*Y$. We can write locally 
\be
\rho^*H =\rho_id\rho^i -\cH\circ \rho dt.
\ee
Then the corresponding coordinate relations read
\be
\dr_iS=\rho _j\dr_i\rho^j-p_i, \quad \dr^iS=\rho_j\dr^i\rho^j,\quad
 \cH'-\cH=\rho_i\dr_t \rho^i -\dr_tS.
\ee
Taken on the graph 
\be
\Delta_\rho=\{(q,\rho(q))\in V^*Y\op\times V^*Y\}
\ee
of the canonical automorphism, the function $S$ plays
the role of a local generating function. For instance, if the graph
$\Delta_\rho$ is coordinatized by $(t,y^i,{y'}^i)$, we obtain the familiar
expression
\be
\cH'-\cH=\dr_tS(t,y^i,{y'}^i).
\ee
\bigskip

\noindent 
{\bf 5. REFERENCE FRAMES}
\bigskip

Every connection $\G$ on the bundle $\pi:Y\to \R$ 
defines a horizontal foliation on
$Y\to \R$ whose leaves are the integral curves of the nowhere
vanishing vector field (\ref{1005}). 
Conversely, let $Y$ admit a horizontal foliation such that, for each
point $y\in Y$, the leaf of this foliation through $y$ is locally determined 
by a section $s_y$ of $V^*Y\to \R$ through $y$. Then, the map
\be
\G:Y\to J^1Y, \qquad \G(y)=j^1_ts_y, \qquad \pi(y)=t, 
\ee
is well defined. This is a connection on $Y\to\R$. 

Given a horizontal foliation on $Y$, there exists the associated atlas of
constant local trivializations of $Y$ such that every leaf of this
foliation is locally generated by the equations $y^i=$const., and the
transition functions
$y^i\to {y'}^i(y^j)$ are independent of the coordinate $t$.$^{16,17}$ 
Two such atlases are said to
be equivalent if their union is also an atlas of constant local
trivializations. They are associated with the same horizontal foliation.
Thus, we have proved the following assertion.

\begin{prop} There is one-to-one correspondence between 
the connections $\G$ on
$Y\to \R$ and the equivalence classes 
of atlases of constant
local trivializations of
$Y$ such that $\G^i=0$ relative to the associated
coordinates, called adapted to $\G$.
\end{prop}

\begin{prop}\label{complcon}
Every trivialization of  $Y\to \R$ yields a complete 
connection on this bundle.
Conversely, every complete connection on $Y\to
\R$ defines a trivialization $Y\cong \R\times M$ such that the
associated coordinates are adapted to $\G$.
\end{prop}

\begin{proof}
Every trivialization of $Y\to\R$ defines the horizontal lift 
$\G=\dr_t$ onto $Y$ of the standard field $\dr_t$ on $\R$ which is
obviously a complete connection on $Y\to\R$.
Conversely, let $\G$ be a complete connection on $Y\to\R$. This
is the generator of the 1-parameter group $G_\G$ which
acts freely on $Y$. The orbits of this action are of course the integral
curves of $\G$. Hence, we obtain a projection 
\be
\pi_\G:Y\to Y/G_\G= M. 
\ee
This projection
together with $\pi:Y\to\R$ defines a trivialization of $Y$.
\end{proof}

One can say that a connection $\G$ on an event space $Y\to\R$ 
describes a reference frame in time-dependent mechanics.
Given a reference frame $\G$, we have the corresponding covariant
differential 
\be
D_\G:J^1Y \to VY,\qquad (t,y^i,y^i_t)\mapsto (t,y^i, \dot y^i=
y^i_t-\G^i).
\ee
Let $s$ be a (local) section of $Y\to\R$.
One can think of $D_\G\circ J^1s$ as being the
relative velocity of the motion $s$ with respect to the reference frame
$\G$. Indeed, 
$D_\G\circ J^1s$ vanishes  identically iff $s$ is
an integral curve of $\G$. 

Let us consider the Hamilton evolution equation (\ref{m96}). For any
connection $\G$ in the splitting (\ref{m96}), there exist holonomic canonical
transformations of $V^*Y$ to the coordinates adapted to $\G$ which 
bring (\ref{m96}) into the familiar Poisson
bracket form
\be
d_{Ht}f= \dr_tf +\{\wt\cH,f\}_V. 
\ee
\bigskip

\noindent 
{\bf 6. LAGRANGIAN POISSON STRUCURE}
\bigskip

In contrast with the Legendre bundle $V^*Y$, the configuration space $J^1Y$
of time-dependent mechanics does not possess any canonical Poisson 
structure in general. 
A Poisson structure on $J^1Y$ depends on the choice of a
Lagrangian
\beq
 L= \cL dt,  \qquad\cL: J^1Y\to\R. \label{a1.81}
\eeq
We will use the notation
$\pi_i=\dr_i^t\cL, \qquad \pi_{ij}=\dr^t_i\dr_j^t\cL.$

Every Lagrangian $L$ (\ref{a1.81}) defines the Legendre map
\beq
\wh L:J^1Y\to V^*Y,\qquad p_i \circ\wh L = \pi_i.\label{a303}
\eeq
The pull-back on $J^1Y$ of the canonical 3-form $\bom$ (\ref{m66}) 
by the Legendre map $\wh L$ (\ref{a303}) reads
\be
\bom_L=\wh L^*\bom =d\pi_i\w dy^i\w dt.
\ee
By means of $\bom_L$, every vertical vector field 
$\vt=\vt^i\dr_i +\dot\vt^i\dr_i^t$ on
$J^1Y\to \R$ yields the 2-form
\be
\vt\rfloor\bom_L= \{[\dot\vt^j \pi_{ji} +\vt^j(\dr_j\pi_i -\dr_i\pi_j)]dy^i-
\vt^i \pi_{ji}dy^j_t\}\w dt.
\ee 
This is one-to-one correspondence, if
 the Lagrangian $\cL$ is regular. 
Indeed, given any 2-form $\f=(\f_idy^i +\dot\f_idy^i_t)\w dt$
on $J^1Y$,
the algebraic equations 
\be
 \dot\vt^j \pi_{ji} +\vt^j(\dr_j\pi_i -\dr_i\pi_j)=\f_i, \qquad 
-\vt^i \pi_{ji}=\dot\f_j
\ee
have a unique solution
\be
\vt^i=-(\pi^{-1})^{ij}\dot\f_j, \qquad \dot\vt^j=(\pi^{-1})^{ji} [\f_i +
(\pi^{-1})^{kn}\dot\f_n(\dr_k\pi_i -\dr_i\pi_k)].
\ee 
In particular, every function $f$ 
on $J^1Y$ determines a vertical vector field
\beq
\vt_f=-(\pi^{-1})^{ij}\dr_j^tf\dr_i + (\pi^{-1})^{ji} [\dr_if +
(\pi^{-1})^{kn}\dr_n^tf(\dr_k\pi_i -\dr_i\pi_k)]\dr_j^t \label{m125}
\eeq
on $J^1Y\to \R$ in accordance with the relation 
\be
\vt_f\rfloor\bom_L = df\w dt.
\ee
Then the Poisson bracket 
\beq
\{f,g\}_Ldt=\vt_g\rfloor\vt_f\rfloor\bom_L, \qquad f,g\in C^\infty(J^1Y),
\label{m121}
\eeq
can be defined on functions on $J^1Y$, and reads
\be
&&\{f,g\}_L= [(\pi^{-1})^{ij} +(\dr_n\pi_k
-\dr_k\pi_n)(\pi^{-1})^{ki}(\pi^{-1})^{nj}](\dr_i^tf\dr_jg-\dr_i^tg\dr_jf) +
 \\
&&\qquad (\dr_n\pi_k
-\dr_k\pi_n)(\pi^{-1})^{ki}(\pi^{-1})^{nj}\dr_i^tf\dr_j^tg. 
\ee
The vertical vector field 
$\vt_f$ (\ref{m125}) is the Hamiltonian vector field
of the function $f$ with respect to the Poisson structure (\ref{m121}).

In particular, if the Lagrangian $\cL$ is hyperregular, that is, the
Legendre map $\wh L$ is a diffeomorphism, the Poisson structure
(\ref{m121}) is obviously isomorphic to the Poisson structure (\ref{m72}) on
the  phase space $V^*Y$. 

The Poisson structure (\ref{m121}) defines the corresponding symplectic
foliation on $J^1Y$ which coincides with the fibration $J^1Y\to \R$. The
symplectic form on the leaf $J^1_tY$ of this foliation is 
$\Om_t= d\pi_i\w dy^i$.$^{18}$

We will see below that the Lagrangian counterpart of Hamiltonian forms
is the Poincar\'e-Cartan form
\be
H_L=\pi_i dy^i -(\pi_iy^i_t-\cL)dt. 
\ee
This is the unique Lepagian equivalent of a Lagrangian $L$ which 
participate in the first variational formula.
Let  
\beq
u=u^t\dr_t +u^i\dr_i, \qquad u^t=0,1, \label{m223}
\eeq
be a vector field on $ Y\to\R$. The first variational formula 
provides the canonical decomposition of the Lie
derivative 
\ben
&&\bL_{J^1 u}L =(J^1u\rfloor d\cL) dt=
(u^t\dr_t +u^i\dr_i +d_tu^i\dr_i^t)\cL dt, \label{1004}\\
&& d_t=\dr_t+y^i_t\dr_i +\cdots, \nonumber
\een
in accordance with the variational task.$^{16,19}$  We have
\beq
J^1 u\rfloor d\cL= 
(u^i-u^ty^i_t)\cE_i +d_t(u\rfloor H_L) \label{m218}
\eeq
where 
\beq
\cE_L=(\dr_i-d_t\dr^t_i)\cL dy^i\w dt\label{983}
\eeq
is the Euler--Lagrange operator for $L$. The kernel $\Ker\cE_L\subset
J^2Y$ of the Euler--Lagrange operator defines
the second order Lagrange equations on $Y$
\beq
(\dr_i-d_t \dr^t_i)\cL=0. \label{b327}
\eeq

\begin{defi}
A connection  $\xi=\dr_t + \xi^i\dr_i + \xi^i_t \dr_i^t$
on the bundle $J^1Y\to \R$
is said to be a 
Lagrangian connection for the Lagrangian
$L$ if it obeys the condition
\beq
\xi\rfloor\Om_L = dH_L \label{2234}
\eeq
which takes the coordinate form
\be
&& (\xi^i - y_t^i)\dr_j\pi_i =0,  \\
&& \dr_iL - \dr_t\pi_i -
\xi^j\dr_j\pi_i - \xi_t^j\dr_j^t\pi_i^t + (\xi^j
- y_t^j)\dr_i\pi_j =0
\ee
relative to the adapted coordinates $(t,y^i,y^i_t,\wh y^i_t,y_{tt})$ on
$J^1J^1Y$.
\end{defi}

In order to clarify the meaning of (\ref{2234}), let us consider the 
Lagrangian 
\be
\ol\cL = \cL + (\wh y_t^i - y_t^i)\pi_i
\ee
on the repeated jet manifold $J^1J^1Y$. The
corresponding Euler--Lagrange operator, called the Euler--Lagrange--Cartan
one, reads
\ben
&& \cE_{\ol L} = [(\dr_i L - \wh d_t\pi_i 
+ \dr_i\pi_j(\wh y_t^j - y_t^j))dy^i + \dr_i^t\pi_j(\wh
y_t^j - y_t^j) dy_t^i]\w dt, \label{2237} \\
&&\wh d_t=\dr_t +\wh y^i_t\dr_i +y^i_{tt}\dr_i^t. \nonumber
\een
Then the condition
(\ref{2234}) is equivalent to the one
$\im \xi\subset \Ker \cE_{\ol L}$, and leads 
 to the first order differential equations
on the jet manifold
$J^1Y$, called the Cartan equations,
\beq
\dr_i^t\pi_j(\wh y_t^j - y_t^j)=0, \qquad
\dr_i L - d_t\pi_i 
+ (\wh y_t^j - y_t^j)\dr_i\pi_j=0. \label{b336}
\eeq
Integral curves of  Lagrangian connections
$\xi$ for $L$ provides classical solution 
$\ol s:\R\to J^1Y$ of these equations.

The restriction of $\cE_{\ol L}$ to the holonomic jet manifold
$J^2Y$ defines the first order Euler--Lagrange operator whose kernel 
is the system of first order Lagrange equations 
\beq
\wh y^i_\la - y^i_\la=0, \qquad (\dr_i-d_\la\dr_i^\la)\cL=0.\label{306}
\eeq
These are equivalent to the second order Lagrange equations
(\ref{b327}), and represent their familiar first order
reduction.

It is easily seen that the first order Lagrange equations 
(\ref{306}) (and consequently the second order
ones (\ref{b327})) are equivalent to the Cartan equations
(\ref{b336}) on the integrable sections $\ol s=J^1s$ of
$J^1Y\to \R$. They are completely equivalent to the Cartan equations in
the case of regular Lagrangians.
\bigskip

\noindent 
{\bf 7. DEGENERATE SYSTEMS}
\bigskip

In time-dependent mechanics, a dynamic equation on the configuration space 
$J^1Y$ is defined to be
a holonomic ($J^2Y$-valued) connection 
\beq
\xi=\dr_t + y^i_t\dr_i + \xi^i_t \dr_i^t \label{a1.30}
\eeq
on the bundle $J^1Y\to \R$. It yields the second order
differential equation on $Y$
\beq
y^i_tt=\xi^i. \label{982}
\eeq

If $\xi$ (\ref{a1.30}) is a
Lagrangian connection for a Lagrangian $L$,
solutions of the dynamic equation (\ref{982})
also satisfy the Lagrange equations (\ref{b327}). 
This is the well-known inverse problem.
If a Lagrangian $L$ is regular, there exists a unique holonomic 
Lagrangian connection 
for $L$. In general, a solution of Lagrange equations is not
necessarily extended to a holonomic Lagrangian connection
and, consequently, is not a
solution of any dynamic equation.  

Turn now to Hamilton equations. 
Every Hamiltonian form $H$ on the Legendre bundle $V^*Y$ defines the
Hamiltonian map 
\be
\wh H: V^*Y\to J^1Y, \qquad y^i_t\circ\wh H=\dr^i\cH.
\ee
Its jet prolongation reads
\be
J^1\wh H: J^1YV^*Y\to J^1J^1Y, \qquad (y^i_t,\wh y^i_t,y^i_{tt})\circ 
J^1\wh H=(\dr^i\cH,y^i_t,d_t\dr^i\cH).
\ee
Given the Hamiltonian connection $\g_H$
(\ref{m57}) for $H$, let consider the composition of morphisms
\beq
J^1\wh H\circ\g_H: V^*Y\to J^2Y,\qquad
 (y^i_t,\wh y^i_t,y^i_{tt})\circ 
J^1\wh H\circ \g_H=(\dr^i\cH,\dr^i\cH,d_{Ht}\dr^i\cH). \label{1002}
\eeq
If the Hamiltonian map $\wh H$ is a diffeomorphism, then 
$J^1\wh H\circ\g_H\circ \wh H^{-1}$ is a dynamic equation.

Let us consider more general condition for solutions of Hamilton
equations to be solutions of the Lagrange and dynamic ones.

Following the general
polysymplectic scheme, we say that a Hamiltonian form
$H$ on $V^*Y$ is associated with a Lagrangian $L$ 
if $H$ obeys the conditions$^{14-16}$
\bea 
&&\wh L\circ\wh H\circ \wh L=\wh L,  \label{m191a}\\
&& p_i\dr^i\cH - \cH=\cL\circ\wh H. \label{m191b}
\eea
It follows from the condition (\ref{m191a})
that $\wh L\circ\wh H$ is the projection operator to 
$Q=\wh L(J^1Y)\subset V^*Y$, called the
Lagrangian constraint space, and $\wh H\circ \wh L$ is the projection
operator to $H(Q)\subset J^1Y$.

If a Lagrangian $L$ is hyperregular, there exists a unique
Hamiltonian form associated with $L$.

Let a Lagrangian $L$ be semiregular, i.e., the pre-image $\wh
L^{-1}(p)$ of any point $p\in Q$ is a
connected submanifold of $J^1Y$. The following assertions 
issue from the corresponding theorems of polysymplectic
formalism.$^{14,16,20}$

\begin{prop} \label{3.22}
All Hamiltonian forms $H$
associated with a semiregular Lagrangian $L$ coincide 
on the Lagrangian constraint space $Q$, and the Poincar\'e--Cartan form
$H_L$ is the pull-back 
of any such a Hamiltonian form $H$ by the Legendre map $\wh L$.
\end{prop}

\begin{prop}\label{3.24}
Let a section $r$ of the bundle $V^*Y\to\R$ be a solution of the Hamilton
equations (\ref{m41}) for a Hamiltonian form $H$ associated
with a semiregular Lagrangian density $L$. If $r$ lives in the Lagrangian
constraint space $Q$,
the section $s=\pi_\Pi\circ r$ of the bundle $Y\to\R$
satisfies the Lagrange equations (\ref{b327}), while its jet prolongation
$\ol s=\wh H\circ r=J^1s$ obeys the Cartan equations
(\ref{b336}).  Conversely, let a section $\ol s$ of
$J^1Y\to\R$ be a solution of the Cartan equations
(\ref{b336}) for a semiregular Lagrangian $L$.
Let $H$ be a Hamiltonian form associated with $L$ so that the corresponding
Hamiltonian map satisfies the condition
\be
\wh H\circ\wh L\circ\ol s=J^1(\pi^1_0\circ\ol s). 
\ee
Then the section $r=\wh L\circ\ol s$ of $V^*Y\to\R$ is a
solution of the Hamilton equations for $H$.
\end{prop}

To prove this Proposition, one can show that, in the case of a semiregular
Lagrangian $L$, the Euler--Lagrange--Cartan operator (\ref{2237}) is the
pull-back
\beq
\cE_{\ol L}= (J^1\wh L)^*\cE_H \label{1001}
\eeq
of the Hamilton operator for a Hamiltonian form $H$ associated with $L$.
In accordance with the relation (\ref{1001}), if $\g_H$ is a
Hamiltonian connection for $H$, the composition $J^1\wh H\circ\g_H$
(\ref{1002}) takes its values into the kernel of the Euler--Lagrange
operator $\cE_L$. Then the morphism $J^1\wh H\circ\g_H\circ\wh L$ restricted
to $\wh H(Q)$ is a local section on $\wh H(Q)\subset J^1Y$ of the affine
bundle $J^2Y\to J^1Y$. If $\wh H(Q)$ is closed, $J^1\wh H\circ\g_H\circ\wh L$
can be extended to a holonomic connection on $J^1Y$. In this case,
projections $\pi_\Pi\circ r$ of integral curves $r$
of the Hamiltonian connection $\g_H$ are also
solutions of a dynamic equation. In particular, this takes place if a
Lagrangian $L$ is almost regular.
              
\begin{defi}  A semiregular Lagrangian density $L$ is said to be almost
regular if (i) the Lagrangian constraint space $Q\to Y$  is a closed
imbedded subbundle $i_Q:Q\hookrightarrow V^*Y$ of the Legendre bundle
$V^*Y\to Y$ and (ii) the Legendre map $\wh L:J^1Y\to Q$ is a bundle.
\end{defi}
    
Since, by Proposition \ref{3.24},
solutions of the Lagrange equations for a degenerate Lagrangian may
correspond to solutions
of different Hamilton equations, we can conclude 
that, roughly speaking, the Hamilton equations involve some additional
conditions in comparison with the Lagrange ones.
Therefore, let us separate a part of the Hamilton equations which are
defined on the Lagrangian constraint space $Q$ in the case of almost
regular Lagrangians.  

Let $H_Q=i_Q^*H$
be the restriction of a Hamiltonian form $H$ 
associated with $L$ to the constraint space $Q$. By virtue of Proposition
\ref{3.22}, this restriction,
called the constrained Hamiltonian form,
is uniquely defined, and $H_L=\wh L^*H_Q$.
For
sections
$r$ of the bundle $Q\to \R$, we can write the constrained Hamilton
equations
\beq
r^*(u_Q\rfloor dH_Q) =0 \label{N44}
\eeq
where $u_Q$ is an arbitrary vertical vector field on
$Q\to \R$.$^{14,16,21}$
In brief, we can identify a vertical vector field $u_Q$ on
$Q\to Y$ with its image $Ti_Q(u_Q)$ and can bring the constrained Hamilton
equations (\ref{N44}) into the form
\beq
r^*(u_Q\rfloor dH) =0 \label{N44'}
\eeq
where $r$ is a section of $Q\to X$ and $u_Q$ is an arbitrary vertical vector
field on $Q\to \R$.
These equations fail to be equivalent to the Hamilton
equations restricted to the constraint space $Q$.

The following two assertions together with Proposition \ref{3.24} give the
relations between Cartan, Hamilton and constrained Hamilton equations
when a Lagrangian is almost regular.$^{16}$

\begin{prop}\label{3.00} For any Hamiltonian form $H$ associated
with an almost regular Lagrangian $L$, every solution $r$ of
the Hamilton equations which lives in the Lagrangian constraint space
is a solution of the constrained Hamilton equations (\ref{N44'}).
\end{prop}

\begin{prop}\label{3.01} A section $\ol s$ of $J^1Y\to X$ is a solution of
the Cartan equations (\ref{b336})
iff $\wh L\circ \ol s$ is a
solution of  the constrained Hamilton equations (\ref{N44'}).
\end{prop}

\begin{rem}
Given a Hamiltonian form $H$ (\ref{b4210}) on $V^*Y$, let us consider the
Lagrangian
\beq
L_H=(p_iy^i_t- \cH)dt \label{1007}
\eeq
on the jet manifold $J^1V^*Y$. It is readily observed that the
Poincar\'e--Cartan form of the Lagrangian $L_H$ coincides with
the Hamiltonian form $H$, and the Euler--Lagrange operator for $L_H$ is
presicely the Hamilton operator for $H$. As a consequence, the Lagrange
equations for $L_H$ are equivalent to the Hamilton equations for $H$.
\end{rem}
 
In the spirit of 
well-known Gotay's algorithm for analyzing constrained systems in symplectic
mechanics,$^{22,23}$
the Lagrangian constraint space $Q$ plays the role of the
primary constraint space. However, one has to apply this algorithm
to each Hamiltonian form $H$ weakly associated with a degenerate Lagrangian 
$L$. If $L$ is semiregular, all these Hamiltonian forms coincide 
on $Q$, but not the corresponding Hamiltonian connections
$\g_H$ (\ref{m57}). 
The necessary condition for a local solution of the Hamilton
equations for a Hamiltonian form
$H$ to live in the Lagrangian constraint space $Q$ is that the 
Hamiltonian connection $\g_H$ is tangent to $Q$ at some
point of $Q$. Given a Hamiltonian form $H$ associated with $L$, we
can express this condition in the explicit form
\bea
&& p_i=\dr^t_i \cL(t, y^j, \dr^j\cH), \label{b4301a}\\
&& (\dr_t +\dr^i\cH\dr_i -\dr_i\cH\dr^i)\rfloor d( p_i-\dr^t_i \cL(t, y^j,
\dr^j\cH))=0. \label{b4301b}
\eea
The equation (\ref{b4301a}) is the coordinate
expression of the relation (\ref{m191a}), and can be taken as 
the equation of the Lagrangian constraint space $Q$. The equation
(\ref{b4301b}) requires that the vector field
$\tau_H$ is tangent to $Q$ at a point with coordinates
$(t,y^i,p_i)$.

In particular, one can apply the description of the quadratic Hamiltonian
systems in polysymplectic formalism$^{16-16}$
to those in time-dependent mechanics. Note that,
since Hamiltonians in time-dependent mechanics are not functions on a phase
space, we cannot apply to them the well-known analysis of the normal
forms$^{24}$
(e.g., quadratic Hamiltonians$^2$) in symplectic
mechanics. 
\bigskip

\noindent 
{\bf 8. CONSERVATION LAWS}
\bigskip

In autonomous mechanics, an integral of motion, by definition, is a function
on the phase space whose Poisson bracket with a
Hamiltonian is equal to zero. This notion cannot be extended 
to  time-dependent mechanics because the Hamiltonian evolution equation
(\ref{m96}) is not reduced to the Poisson bracket. 

We start from conservation laws in Lagrangian mechanics.
To obtain differential
conservation laws, we use the first
variational formula (\ref{m218}).$^{16,19}$
On-shell, this 
leads to the weak identity
\beq
J^1 u\rfloor d\cL\ap -d_t\cT \label{m219}
\eeq
where
\beq
\cT=\pi_i(u^ty^i_t -u^i) -u^t\cL \label{m225}
\eeq
is the current along the vector field $u$ (\ref{m223}). 
If the Lie derivative $\bL_{J^1 u}L$ (\ref{1004})
vanishes, we have the conservation law 
\be
0\ap - d_t[\pi_i(u^ty^i_t-u^i )-u^t\cL]. 
\ee
This is brought into the differential conservation law
\be
0 \ap -\frac{d}{dt}(\pi_i\circ s(u^t \dr_t s^i-u^i\circ s) -u^t\cL\circ s) 
\ee
on solutions $s$ of the Lagrange equations. A glance at
this expression shows that, in time-dependent
mechanics, the conserved current (\ref{m225})
plays the role of an integral of motion. 

Every symmetry current (\ref{m225}) along a
vector field $u$ (\ref{m223}) on $Y$ 
can be represented as a superposition 
of the N\"other current along a
vertical vector field $\vt$ and  of the energy current
along some connection $\G$ (\ref{1005}) on $Y\to\R$, where 
$u=\vt+\G$.$^{16,25}$

If $\vt$ is a vertical vector field, the weak identity (\ref{m219}) reads
\be
(\vt^i\dr_i +d_t\vt^i \dr^t_i)\cL \ap d_t(\pi_i\vt^i). 
\ee
If the Lie derivative of $L$ along $\vt$ equals  zero, we have the integral
of motion $\cT=\pi_i\vt^i$.

In the case of a connection
$\G$ (\ref{1005}), the weak identity (\ref{m219}) takes the form
\beq
(\dr_t +\G^i\dr_i +d_t\G^i\dr_i^t)\cL \ap - d_t(\pi_i(y^i_t -\G^i) -\cL),
\label{m227}
\eeq
where one can think of
\beq
\cT_\G= \pi_i(y^i_t -\G^i) -\cL \label{m228}
\eeq
as being the energy function with respect to the reference
frame $\G$. In particular, the energy conservation law (\ref{m227})
written relative to the coordinates adapted to $\G$ takes the familiar form
\beq
\dr_t\cL =- d_t(\pi_iy^i_t -\cL). \label{m229}
\eeq

To discover conservation laws within the framework of  Hamiltonian
formalism, let us
consider the Lagrangian 
(\ref{1007}) on $J^1V^*Y$,  
and apply the first variational formula (\ref{m218}) to it.$^{19}$

Given a vector field (\ref{m223}) on the event bundle $Y$, its canonical
lift  onto $V^*Y$ reads
\be
\wt u=u^t\dr_t + u^i\dr_i - \dr_i u^j p_j\dr^i, \qquad u^t=0,1. 
\ee
Substituting this vector field into the weak identity (\ref{m219}), we obtain
\beq
-u^i\dr_i\cH -u^t\dr_t\cH  
+p_id_t u^i \ap - d_t(-p_iu^i +u^t\cH) \label{m230} 
\eeq
for the current
\beq
\wt\cT= -p_iu^i + u^t\cH. \label{b4306}
\eeq
 
In the case of $u=\G$, 
the weak identity (\ref{m230})  takes the form
\be
-\dr_t\cH -\G^i\dr_i\cH +p_id_t \G^i \ap -d_t\wt\cH_\G
\ee
where $\wt\cH_\G=\cH -p_i\G^i$ is
the Hamiltonian function in the splitting (\ref{m46}).

The following assertion shows that the Hamiltonian function $\wt\cH_\G$
is the Hamiltonian counterpart of the Lagrangian energy function
$\cT_\G$ (\ref{m228}) in the case of semiregular Lagrangians.$^{16}$

\begin{prop} Let a 
Hamiltonian form $H$ on the Legendre bundle $V^*Y$  be associated
with a semiregular Lagrangian  $\cL$ on $J^1Y$. Let $r$ be a solution
of the Hamilton equations (\ref{m41}) for $H$ which lives in
the Lagrangian constraint space $Q$ and $s$ the corresponding
solution  of the
Lagrange equations for $L$. Then, we have 
\be
\wt\cT (r)=\cT( \wh H\circ r),\qquad
\wt\cT (\wh L\circ J^1s) =\cT(s)
\ee
where $\cT$ is the current (\ref{m225}) on $J^1Y$ and $\wt\cT$ is the
current (\ref{b4306}) on $V^*Y$.
\end{prop}

Therefore, we can treat $\wt\cH_\G$ as the energy function
with respect to the reference
 frame $\G$. In particular, if
$\G^i= 0$, we obtain the well-known energy conservation law
\be
\dr_t\cH \ap d_t\cH
\ee
which is the Hamiltonian variant of the Lagrangian one (\ref{m229}). 
\bigskip

\noindent 
{\bf 9. RELATIVISTIC MECHANICS}
\bigskip

Let us consider a mechanic system whose event space $Z$ has no
fibration $Z\to \R$ or admits
different such fibrations. We come to relativistic mechanics where a 
configuration space is the 
jet manifold of 1-dimensional submanifolds of $Y$ that generalizes the
notion of jets of sections of a bundle.$^{16,26,27}$

Let $Z$ be a manifold of dimension $m+n$. 
The 1-order jet manifold $J^1_nZ$ of $n$-dimensional submanifolds of $Z$
comprises the equivalence classes $[S]^1_z$  of
$n$-dimensional imbedded submanifolds of $Z$ which pass through $z\in
Z$  and which are tangent to each other at $z$.
It is provided
with a manifold structure as follows.

Let $Y\to X$ be an $(m+n)$-dimensional bundle over an
$n$-dimensional base $X$ and $\Phi$  an imbedding of $Y$ into $Z$.
Then there is the natural injection
\be
&& J^1\Phi: J^1Y\to J^1_nZ, \\
&& j^k_xs \mapsto [S]^k_{\Phi(s(x))}, \qquad S=\im (\Phi\circ s), 
\ee
where $s$ are sections of $Y\to X$.
This  injection  defines a chart on
$J^1_nZ$.  Such charts cover the set $J^1_nZ$,
and transition functions between these charts are differentiable. They
provide 
$J^1_nZ$ with the structure of a finite-dimensional manifold.

Hereafter, we will use the following coordinate atlases on the
jet manifold $J^1_nZ$ of submanifolds of $Z$.
Let $Z$ be endowed with a manifold atlas with coordinate charts
\beq
(z^A), \qquad A=1,\ldots,n+m.  \label{5.59}
\eeq
Though $J^0_nZ$, by definition, is diffeomorphic to $Z$, let us provide
$J^0_nZ$ with the atlas 
obtained by replacing  every chart $(z^A)$ on a domain
$U\subset Z$ with the
charts on the same domain $U$ which
correspond to the different partitions of the collection $(z^A)$ in
collections of $n$ and $m$ coordinates. We denote these coordinates by
\beq
(x^\la,y^i), \qquad \la=1,\ldots,n,  \qquad i=1,\ldots,m.\label{5.8}
\eeq
The transition functions between the coordinate charts (\ref{5.8}) of
$J^0_nZ$ associated with the coordinate chart
(\ref{5.59}) of $Z$ are reduced simply to exchange between 
coordinates $x^\la$ and $y^i$. 
Transition functions between arbitrary coordinate charts of the
manifold $J^0_nZ$  take the form
\beq
\wt x^\la = \wt g^\la (x^\mu, y^j), \qquad \wt y^i = \wt f^i (x^\mu, y^j).
\label{5.26}
\eeq

Given the coordinate atlas (\ref{5.8}) of the manifold $J^0_nZ$, the
jet manifold $J^1_nZ$ of $Z$ is endowed with the adapted
coordinates $(x^\la,y^i_\la)$. 
Using the formal total derivatives $d_\la=\dr_\la +y^i_\la\dr^i +\cdots$,
one can write the
 transformation rules for these coordinates in the
following form.
 Given the coordinate transformations (\ref{5.26}), it is easy to find that
\beq
d_{\wt x^\la} = \left[d_{\wt x^\la}
g^\al (\wt x^\la, \wt y^i)\right]d_{x^\al}. \label{5.35}
\eeq
Then we have
\beq
\wt y^i_\la =\left[\left(\frac{\dr}{\dr\wt x^\la} + \wt y^p_\la
\frac{\dr}{\dr\wt y^p}\right)
g^\al (\wt x^\la, \wt y^i)\right]\left(\frac{\dr}{\dr x^\al} +y^j_\al
\frac{\dr}{\dr y^j}\right)\wt f^i (x^\mu, y^j). \label{5.36}
\eeq

\begin{rem}
Given a manifold $Z$, there is one-to-one correspondence
between the jets
$[S]^1_z$ at a point $z\in Z$ and the $n$-dimensional vector subspaces
of the
tangent space $T_zZ$:
\be
[S]^1_z \mapsto \dot x^\la(\dr_ \la +y^i_\la([S]^1_z)\dr_i).
\ee
The bundle $J^1_mZ\to Z$ possesses the structure group
$GL(n,m;\R)$ of linear transformations
of the vector space $\R^{m+n}$ which preserve the 
subspace
$\R^n$. Its typical fibre is the Grassmann
manifold $GL(n+m;\R)/GL(n,m;\R)$
of $n$-dimensional vector subspaces
of the vector space $\R^{m+n}$.
In particular, if $n=1$, the fibre coordinates $y^i_0$ of $J^1_1Z\to Z$
with the transition
functions (\ref{5.36}) are exactly the standard
coordinates of the projective space ${\bf RP}^m$.
\end{rem}

When $n=1$, the formalism of jets of submanifolds
provides the adequate mathematical description of relativistic mechanics
as follows.

Let $Z$ be a $(m+1)$-dimensional manifold
equipped with an atlas of coordinates $(z^0,z^i)$, $i=1,\ldots,m,$
(\ref{5.8}) with the transition functions (\ref{5.26}) which take the form
\beq
z^0\to \wt z^0(z^0, z^j), \qquad z^i\to \wt z^i(z^0, z^j). \label{b5.1}
\eeq
The coordinates $z^0$ in different charts of $Z$ play the role of temporal
ones. 

Let $J^1_1Z$ be the jet manifold of 1-dimensional submanifolds of $Z$. This
is provided with the adapted coordinates
$(z^0,z^i,z^i_0)$. Then one can think of $z^i_0$
as being the coordinates of non-relativistic
velocities. Their transition
functions are obtained as follows.

Given the coordinate transformations (\ref{b5.1}), the total derivative
(\ref{5.35}) reads
\be
d_{\wt z^0}= d_{\wt z^0}(z^0)d_{z^0}=\left(\frac{\dr z^0}{\dr \wt z^0} +
\wt z^k_0\frac{\dr z^0}{\dr \wt z^k}\right)d_{z^0}.
\ee
In accordance with the relation (\ref{5.36}), we have
\be
\wt z^i_0=  d_{\wt z^0}(z^0)d_{z^0}(\wt z^i) =
\left(\frac{\dr z^0}{\dr \wt z^0} +\wt z^k_0\frac{\dr z^0}{\dr \wt
z^k}\right) \left(\frac{\dr \wt z^i}{\dr z^0} +
z^j_0\frac{\dr \wt z^i}{\dr z^j}\right).
\ee
The solution of this equation is
\be
\wt z^i_0= \left(\frac{\dr \wt z^i}{\dr z^0} +
 z^j_0\frac{\dr \wt z^i}{\dr z^j}\right)\mbox{\large /}
 \left(\frac{\dr \wt z^0}{\dr z^0} + z^k_0\frac{\dr \wt z^0}{\dr z^k}\right).
\ee
This is the transformation law of non-relativistic
velocities, which illustrates that the
jet bundle $J^1_1Z\to Z$ is not affine, but projective.

To obtain the relation between non-relativistic and relativistic velocities,
let us consider the tangent bundle $TZ$ equipped with
the induced coordinates $(z^0,z^i,\dot z^0, \dot z^i)$.
There is the morphism 
\beq
\rho: TZ\op\to J^1_1Z, \qquad
z^i_0\circ \rho = \dot z^i/\dot z^0. \label{b5.3}
\eeq
It is readily observed that the coordinate transformation laws of $z^i_0$ and
$\dot z^i/\dot z^0$ are the same.
Thus, one can think of the
coordinates $(\dot z^0, \dot z^i)$ as being relativistic
velocities.

\begin{rem}
Note that the similar morphism $\R^{m+1}\to {\bf
RP}^{m+1}$
provides the projective space ${\bf RP}^4$ with the standard coordinate
charts.
\end{rem}

The morphism (\ref{b5.3}) is a surjection. Let us assume that the tangent
bundle is equipped with a pseudo-Riemannian metric $g$ and $Q_z\subset T_zZ$
is the hyperboloid given by the relation
\be
g_{\m\nu}(z)\dot z^\m\dot z^\nu = 1, \qquad \m,\nu=0,1,\cdots m.
\ee
The union of these hyperboloids over $Z$
\be
Q=\op\cup_{z\in Z}Q_z =Q^+\cup Q^-
\ee
is the union of two connected imbedded subbundles of $TZ$.
Then the restriction
of the morphism (\ref{b5.3}) to each of this subbundle is an injection
of $Q$ into $J^1_1Z$.

Let us consider the image of this injection in the fibre of $J^1_1Z$
over a point $z\in Z$. There are coordinates $(z^0,z^i)$ in a neighbourhood
around $z$ such that the pseudo-Riemannian metric $g(z)$ at $z$ comes
to the pseudo-Euclidean one $g(z)={\rm diag}(1,-1,\cdots, -1)$.
In this coordinates the 
hyperboloid $Q_z\subset T_zZ$ is given by the relation
\be
(\dot z^0)^2 -\op\sum_i (\dot z^i)^2=1.
\ee
This is the union of the subsets $Q^+_z$ where $z^0>0$ and $Q^-_z$ where
$z^0<0$. The image $\rho(Q^+_z)$ is given by the coordinate relation
\be
\op\sum_i (z^i_0)<1.
\ee
>From the physical viewpoint, this relation means that non-relativistic
velocities are bounded
in accordance with Special Relativity.

\end{document}